\documentclass[12pt]{article}

\usepackage[utf8]{inputenc}
\usepackage{fullpage}
\usepackage[small,bf]{caption}
\usepackage{todonotes}

\usepackage{float}
\usepackage{geometry}
\usepackage{verbatim}
\usepackage{graphicx}
\usepackage{enumerate, enumitem}
\usepackage{booktabs}
\usepackage[
colorlinks=true,
citecolor=blue,
linkcolor=magenta,
urlcolor=cyan]{hyperref}
\usepackage{url}
\usepackage{xcolor}
\usepackage{tikz}
\usepackage{xparse}
\usepackage{amssymb,amsthm,amsmath,amsfonts,amsthm,mathtools}
\newcommand{\eg}{{\it e.g.}}
\newcommand{\ie}{{\it i.e.}}

\newcommand{\BEQ}{\begin{equation}}
\newcommand{\EEQ}{\end{equation}}
\newcommand{\BEAS}{\begin{eqnarray*}}
\newcommand{\EEAS}{\end{eqnarray*}}

\newcommand{\reals}{{\mbox{\bf R}}}

\newcounter{algorithmctr}
\renewcommand{\thealgorithmctr}{\arabic{algorithmctr}}
   {\mbox{}\\*[\parskip]\begin{minipage}{\linewidth}%
       \refstepcounter{algorithmctr}\begin{list}{}{%
       \setlength{\rightmargin}{0\linewidth}%
       \setlength{\leftmargin}{.05\linewidth}}%
       \rmfamily\small
       \item[]{\setlength{\parskip}{0ex}\hrulefill\par%
        \nopagebreak{\bfseries\textsf{Algorithm \thealgorithmctr~}}}}%
   {{\setlength{\parskip}{-1ex}\nopagebreak\par\hrulefill\\*[2ex]\par}%
   \end{list}\end{minipage}}

\DeclarePairedDelimiter{\normabs}{\lvert}{\rvert}

\makeatletter
\def\abs{\@ifstar{\normabs}{\normabs*}}
\makeatother

\usetikzlibrary{math}

\newcommand{\tr}{\mathop{\bf tr}}

\usepackage{textcomp}

\begin{document}

\title{Computing Tighter Bounds on the $n$-Queens Constant via Newton's Method}
\author{Parth Nobel\thanks{Corresponding author: \texttt{ptnobel@stanford.edu}} \and Akshay Agrawal \and Stephen Boyd \and \\
Department of Electrical Engineering, Stanford University \\
    350 Jane Stanford Way \\
    Stanford, CA 94305 \\
}
\date{
    December 6, 2021
    }

\maketitle

\begin{abstract}
In recent work Simkin shows that bounds
on an exponent occurring in the famous $n$-queens problem
can be evaluated by solving convex optimization problems,
allowing him to find bounds far tighter than previously known.
In this note we use Simkin's formulation, a sharper bound developed by Knuth,
and a Newton method that scales to large problem instances, to find even
sharper bounds.
\end{abstract}

\section{Introduction}
Let $\mathcal Q(n)$ denote the number of ways that $n$ queens can be arranged on
an $n \times n$ chessboard in such a way that none is threatening
another, \ie, no two queens are in the same row, column, or diagonal.
Recent work by Simkin \cite{simkin2021number} has shown that 
\[
\lim_{n \to \infty} \frac{ \mathcal Q(n)^{1/n}}{n} = e^{-\alpha},
\]
where $\alpha$ is a constant, that we refer to as the $n$-queens constant,
characterized as the optimal value of an infinite
dimensional convex optimization problem.
For background on the problem and previously derived bounds on
$Q(n)$, see \cite{bell2009survey}.

In \cite{simkin2021number}, Simkin establishes that $\alpha \in [1.94, 1.9449]$, a strong
tightening of the best previously known bounds $\alpha \in [1.58,3]$ \cite{luria2017new,luria2021lower}.
His method finds lower and upper bounds by solving two convex optimization problems.
Knuth later formulated another convex optimization problem
which also gives an upper bound on $\alpha$ \cite{knuth2021xqueens}.

In this note we solve the convex optimization problems associated with
Simkin's lower bound and Knuth's upper bound, using a version of Newton's
method that scales to large problem instances, to establish that
\[
\alpha \in [1.944000752, 1.944001082].
\]
This agrees with previous conjectures that 
$\alpha \approx 1.944$ \cite{zhang2009counting}.
In terms of the gap, \ie, difference of known upper and lower bounds, Simkin improved
it from the previous value around $1.4$ to 
around $5\times 10^{-3}$, and we have improved that to around $3.3 \times 10^{-7}$.

Simkin's numerical lower bound is found as a lower bound on the optimal
value of a convex optimization problem whose optimal value is a lower bound on 
the $n$-queens constant.
This problem is parameterized by $n$, the size of the chessboard used to interpret
the problem.
We let $L_n$ denote the optimal objective value of this problem.
In \cite{knuth2021xqueens}, Knuth introduces a convex optimization problem
whose optimal value is an upper bound on the $n$-queens constant.
It is also parameterized by $n$, and we let $U_n$ denote its optimal value.
Simkin's numerical upper bound is found by solving a related convex optimization problem
which upper bounds Knuth's problem.

Simkin's numerical bounds are a lower bound on $L_{17}$ and an upper bound on $U_{12}$,
obtained by approximately solving these two problems.
These problems involve a few hundred variables and constraints. 
In contrast, we use Newton's method to solve the problems, which has two
advantages.  First, we solve the problem to high accuracy, so almost nothing is
lost when we move from an approximate solution to a lower or upper 
bound on the optimal value. Second, our method scales to much larger $n$, which gives
us tighter bounds on the $n$-queens constant.

Our numerical bounds are $L_{2048}$ and $U_{1024}$, obtained by solving 
the two problems to high accuracy.
The lower bound problem contains almost 17 million variables and over 12000 constraints;
the upper bound problem contains over 4 million variables and over 14000 constraints.
In this note we explain how a version of Newton's method can be used
to solve such large problems.
(The numbers $n=1024$ and $n=2048$ are chosen as the largest powers of two
that we can evaluate on the desktop computer we used to carry out the 
computations.)

\section{The convex problems}
The bounds $L_n$ and $U_n$ are the optimal values of 
convex optimization problems of the form
\BEQ\label{e-prob}
\begin{array}{ll} \mbox{minimize} & f(x) \\
\mbox{subject to} & Ax=b,
\end{array}
\EEQ
where $x \in \reals^p$ is the variable, $A \in \reals^{q \times p}$, $b\in \reals^q$
specify the constraints, and 
the objective function $f:\reals_{++}^p \to \reals$ is smooth and strictly convex.
($\reals_{++}$ denotes the set of positive numbers.)
These problems are feasible, and so have a unique solution.

In these optimization problems, $f$, $A$, and $b$ are parametrized by $n$,
but to lighten the notation we suppress this dependence 
in our description of the method.
Full descriptions of $f$, $A$, and $b$ for the lower bound and upper
bound problems are given in the appendix.  
Here, we summarize some of their attributes.

\paragraph{Lower bound problem.}
For the lower bound problem, we have $p = 4n^2 + 4n$ variables and $q = 6n - 1$ constraints.
The objective $f$ is separable, \ie, a sum of functions of $x_i$, so its
Hessian $\nabla^2 f(x)$ is diagonal.
The constraint coefficient matrix $A$ is full rank and sparse, with 
at most $4n$ nonzero entries in each row, and at most $4$ nonzero entries in each column.
The entries of $A$ are all $0$ or $1$.

\paragraph{Upper bound problem.}
For the upper bound problem, we have
$p= 4n^2 + 8n - 4$ variables and $q = 14n - 6$ constraints.
The objective $f$ is block separable, a sum of 
functions of pairs of variables, where the pairs are disjoint,
so its Hessian $\nabla^2 f(x)$ is block diagonal, 
with $1 \times 1$ and $2\times 2$ blocks.
Here too $A$ is full rank and sparse, with at most $4$ nonzero entries per column,
and at most $2n+1$ nonzero entries per row.  Its entries are 
all either $0$ or $1$ or $2n$.

\section{Infeasible start Newton method}
In this section we summarize the infeasible start Newton method
described in \cite[\S10.3.2]{bv2004convex} (which also contains a convergence proof),
and explain how to compute the search directions in a scalable way.
We also discuss how to compute appropriate bounds on optimal values of
the problems.

\subsection{Optimality condition and residuals}
The necessary and sufficient optimality conditions for \eqref{e-prob} are
\[
\nabla f(x) + A^T \nu = 0 , \qquad A x - b = 0,
\]
where $\nu \in \reals^q$ is a dual variable or Lagrange multiplier.
For $x \in \reals_{++}^p$ and $\nu \in \reals^q$ 
we define the dual and primal residuals as
\BEQ 
r_d(x, \nu) = \nabla f(x) + A^T \nu, \qquad
r_p(x, \nu) = A x - b,
\EEQ
and the (primal-dual) residual $r(x, \nu) = (r_d(x, \nu), r_p(x, \nu))$.
Thus the optimality condition can be expressed as $r(x,\nu)=0$.

\subsection{Infeasible start Newton method}
The method is iterative, with iterates denoted as $(x^{(k)}, \nu^{(k)})$,
where $k$ is the iteration number.
The iterates will satisfy $x^{(k)} \in \reals_{++}^p$, so the 
residual $r^{(k)} = r(x^{(k)},\nu^{(k)})$
is defined (and will converge to zero as $k \to \infty$).
We initialize our algorithm with
$x^{(0)} \in \reals_{++}^p$, which need not satisfy $A x^{(0)} = b$.

\paragraph{Newton step.}
For the $k$th iterate the Newton step $(\Delta x^{(k)},\Delta \nu^{(k)})\in
\reals^p \times \reals^q$ is the solution of the linear equations
\[
r(x^{(k)}, \nu^{(k)}) + Dr (x^{(k)}, \nu^{(k)}) (\Delta x^{(k)},\Delta \nu^{(k)}) = 0,
\]
where $Dr$ is the derivative or Jacobian of the residual.
(We will show later that these equations always have a unique solution.)
The lefthand side is the first order Taylor approximation of
$r(x^{(k)}+\Delta x^{(k)},\nu^{(k)}+ \Delta \nu^{(k)})$, so
if the Newton step is added to the current iterate, we 
obtain primal and dual variables for which the Taylor approximation is zero.

We write the equations defining the Newton step as
\BEQ\label{e-step-matrix-form}
\left[ \begin{array}{cc} 
\nabla^2 f(x^{(k)}) & A^T \\
A & 0
\end{array}\right]
\left[ \begin{array}{c}
\Delta x^{(k)} \\ \Delta \nu^{(k)}
\end{array}\right] =
- \left[ \begin{array}{c}
r_d (x^{(k)},\nu^{(k)}) \\
r_p (x^{(k)},\nu^{(k)})
\end{array}\right].
\EEQ
The coefficient matrix is invertible, since its top left block is invertible and
its bottom left block is wide and full rank; see, \eg,
\cite[\S16.2]{bv2018vmls} or \cite[\S10.1.1]{bv2004convex}.

From the top block equations in (\ref{e-step-matrix-form}) we have 
\BEQ\label{e-dxk}
\Delta x^{(k)} = - \nabla^2 f(x^{(k)})^{-1}\left( 
r_d (x^{(k)},\nu^{(k)}) + A^T \Delta \nu^{(k)} \right).
\EEQ
Substituting this into the bottom block of equations we obtain the set 
of equations
\BEQ\label{e-delta-nu}
    \left( A \nabla^2 f(x^{(k)})^{-1} A^T\right)
\Delta \nu^{(k)} = r_p(x^{(k)}, \nu^{(k)}) - 
    A \nabla^2 f(x^{(k)})^{-1} r_d(x^{(k)}, \nu^{(k)}),
\EEQ
with positive definite coefficient matrix $S =  A \nabla^2 f(x^{(k)})^{-1} A^T$.
To find the Newton step, we first solve the set of equations \eqref{e-delta-nu}
to obtain $\Delta \nu^{(k)}$, and then evaluate $\Delta x^{(k)}$ using
\eqref{e-dxk}.

\paragraph{Line search and update.}
The next iterate has the form
\[
x^{(k+1)} = x^{(k)} + t^{(k)} \Delta x^{(k)}, \qquad
\nu^{(k+1)} = \nu^{(k)} + t^{(k)} \Delta \nu^{(k)},
\]
where $t^{(k)}$ is a positive step length.
Choosing $t^{(k)}$ is referred to as a line search.
Our line search is one specifically for the infeasible start 
Newton method, described in \cite[\S 10.3.2]{bv2004convex}; for
more general discussion of line search methods, see, \eg,
\cite[Chap.~3]{nw2006numopt}.

To find $t^{(k)}$ we first find $\tilde t= \min(0.95t_\text{max}, 1)$, where 
\[
t_\text{max} = \min \left. \left\{ \frac{x^{(k)}_i}{-\Delta x^{(k)}_i} \; \right| \;
\Delta x^{(k)}_i<0 \right\}
\]
is the largest possible step
for which $x^{(k)}+t \Delta x^{(k)}  \in \reals_{+}^p$.
We take $t^{(k)} = \beta^\ell \tilde t$, where $\beta \in (0,1)$ is a parameter
and $\ell$ is the smallest positive integer for which 
\[
\| r(x^{(k)} + \beta^\ell \tilde t \Delta x^{(k)}, \nu^{(k)}+ \beta^\ell \tilde t
\Delta \nu^{(k)}) \|_2 \leq
(1-\alpha \beta^\ell \tilde t)
\| r(x^{(k)}, \nu^{(k)}) \|_2
\]
holds, where $\alpha \in (0,1/2)$ is a parameter.
(It can be shown that such an integer exists \cite[\S10.3.1]{bv2004convex}.)
If $\|r(x^{(k+1)}, \nu^{(k+1)})\|_2 < \epsilon$, we terminate, where $\epsilon$ is a positive
tolerance. 

We use the common line search parameter values 
$\alpha = 0.01$ and $\beta = 0.9$, and the tolerance $\epsilon = 10^{-9}$, 
which is far smaller than would be needed in any engineering or statistics application.

\paragraph{Efficient computation.}
The computational effort in each step is, predominantly, solving the
linear equations \eqref{e-delta-nu}.
Since we intend to use our algorithm on problem instances where forming
and storing the $q \times q$ 
matrix $S = A \nabla^2 f(x^{(k)})^{-1} A^T$ is not practical,
we use an indirect iterative method to solve these equations \cite{saad2003iterative}.
There are many such methods, mostly based on Krylov subspaces, such as
conjugate gradients \cite{hestenes1952cg}.
The particular method we use is MINRES \cite{paige1975solution}.
Like other indirect methods, it requires only a method to 
evaluate the mapping $y \mapsto Sy$ for a vector $y$.  
We evaluate this mapping as
\[
Sy = A \left( \nabla^2 f(x^{(k)})^{-1} \left( A^Ty \right)\right),
\]
\ie, successive multiplications by $A^T$, $\nabla^2 f(x^{(k)})^{-1}$, and $A$,
without forming or storing the matrix $S$.
(Simply storing $S$ for the specific problems we will solve 
would require many terabytes of memory.)

\subsection{Bounds on optimal values}
To bound the $n$-queens constant $\alpha$, we need a lower bound on $L_n$ and
an upper bound on $U_n$.
Newton's method is able to solve the lower bound and upper bound problems 
to high accuracy, so the optimal values we compute could simply be rounded down
or up to obtain these bounds.
Here we discuss ways to more carefully compute these bounds on the optimal
values.
To find an upper bound on $U_n$, it suffices to find a feasible $x$ and evaluate 
the objective at that point.

\paragraph{Lower bounds on $L_n$.}
We follow Simkin and use standard Lagrangian duality to find a lower bound on $L_n$.
The dual function of \eqref{e-prob} is
\[
h(\nu) = \nu^T b - f^*(A^T\nu),
\]
where $f^*$ is the conjugate function of $f$ \cite[\S3.3]{bv2004convex}.
For \emph{any} $\nu$, $h(\nu)$ is a lower bound on the optimal value of the problem
\eqref{e-prob}.
For the lower bound problem we can explicitly find $f^*$ as
\[
f^*(y) = \sum_{i=1}^p \exp \left(y_i - 1\right) - 4 \log n - 2 \log 2 - 3
\]
via considering the structure of $f$ as given in \S\ref{s-lower-bound-problem}
and a straightforward application of results given in \cite[\S3.3.1]{bv2004convex}.
To obtain a lower bound on $L_n$, we solve the problem using Newton's method 
and then evaluate $h(\nu)$ for the $\nu$ found.
(Since we solve these problems to high accuracy, the lower bound on $L_n$ obtained
is very close to the upper bound on $L_n$ found, which is $f(x)$.)

\paragraph{Rational approximation.}
The bounds described above are found using floating point computations.
To make the upper bound fully precise, we find a rational approximation of $x$
that is \emph{exactly} feasible and evaluate the objective, carefully using
an upper bound on the (transcendental) objective function.
For the lower bound, we find a rational approximation of $\nu$ and 
evaluate a lower bound on the dual function.  We have not taken these
steps, because our floating point solutions are so accurate that it would 
have a negligible effect on our final numerical bounds.

\section{Results}

\paragraph{Lower bound.}
We computed $L_{2048}$ by solving a problem with $p=16785408$ variables and 
$q=12288$ constraints. This required 21 iterations, with a total time of
around 517 seconds on a M1 Mac Mini.
We obtained the lower bound 
\[
    1.944000752019729 = L_{2048}.
\]

\paragraph{Upper bound.}
We computed $U_{1024}$ by solving 
a problem with $p=4202492$ variables and $q=14332$ constraints.
To speed up finding the solution, we started Newton's method from $(x,\nu)$ that solve
the approximate upper bound problem \eqref{e-approx-upper-bnd}, which has a simpler
objective and so was faster to compute.
This required 6 iterations to solve the approximate problem, and a further 7 iterations
to solve the exact problem.  The total time was around 56 seconds on a M1 Mac Mini.
We obtained the upper bound
\[
    1.9440010813092217 = U_{1024}.
\]

Our code is available at \url{https://github.com/cvxgrp/n-queens}.

\section*{Acknowledgements}
We thank Don Knuth for introducing us to this problem and substantial help
in navigating the various lower and upper bounds.
We thank Michael Simkin for his comments and suggestions on an earlier draft
of this note.

Parth Nobel was supported in part by the National Science Foundation Graduate
Research Fellowship Program under Grant No. DGE-1656518. Any opinions,
findings, and conclusions or recommendations expressed in this material are those of the
author(s) and do not necessarily reflect the views of the National Science Foundation.
This research was partially supported by ACCESS --- AI Chip Center for Emerging
Smart Systems, sponsored by InnoHK funding, Hong Kong SAR.

\clearpage
\bibliographystyle{plain}
\bibliography{n-queens-constant}

\clearpage
\appendix
\section{Details of convex problems}

In this appendix we give the details of the lower and upper bound problems,
as well as an approximate upper bound problem.
We define the variables in their natural notation, leaving 
it to the reader to re-arrange these into a single vector variable $x$.
In a similar way, we describe the linear constraints in their natural notation,
leaving it to the reader to translate these into $Ax=b$.

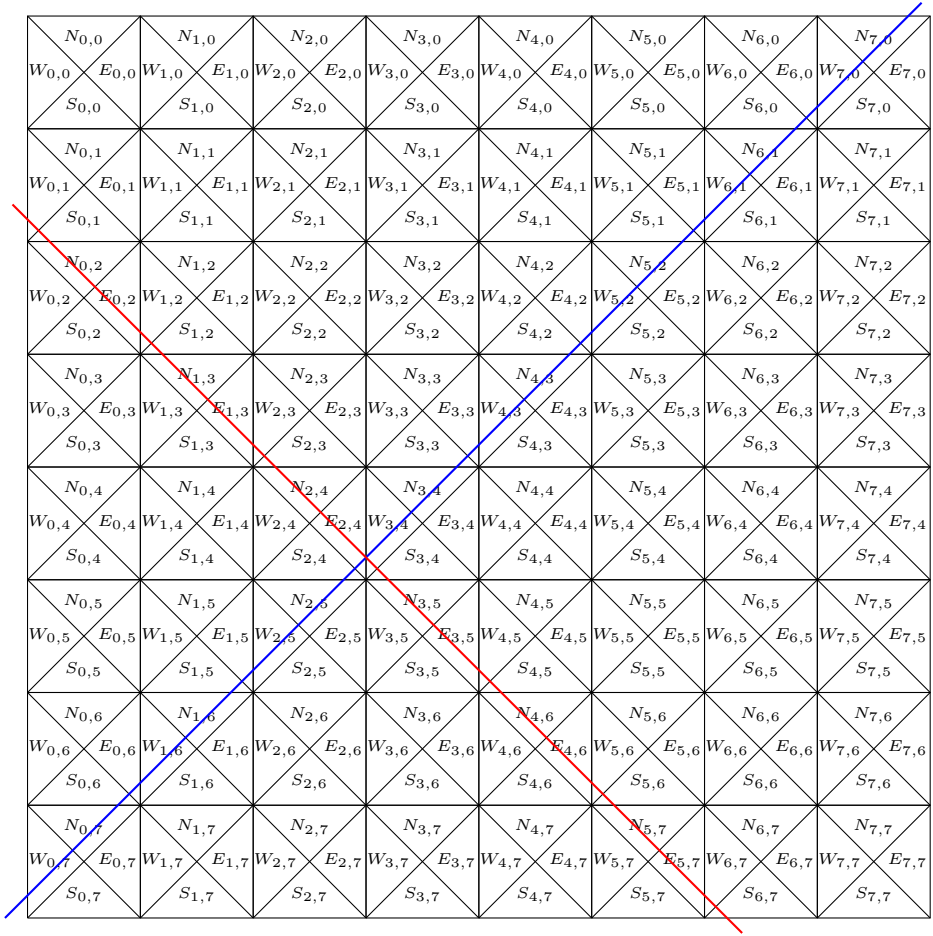
\begin{figure}
    \centering
\begin{tikzpicture}
    \foreach \i in {1,...,8}{
        \foreach \j in {1,...,8}{
            \tikzmath{
                int \iidx, \jidx;
                \iidx=\i-1;
                \jidx=8-\j;
            }
            \node at (1.5*\i+.75, 1.5*\j+1.2) {\tiny \(N_{\iidx, \jidx}\)};
            \node at (1.5*\i+.75, 1.5*\j+.3) {\tiny \(S_{\iidx, \jidx}\)};
            \node at (1.5*\i+.3, 1.5*\j+.75) {\tiny \(W_{\iidx, \jidx}\)};
            \node at (1.5*\i+1.2, 1.5*\j+.75) {\tiny \(E_{\iidx, \jidx}\)};
            \draw (1.5*\i,1.5*\j)
                    -- (1.5*\i, 1.5*\j+1.5)
                    -- (1.5*\i+1.5, 1.5*\j+1.5)
                    -- (1.5*\i+1.5, 1.5*\j)
                    -- (1.5*\i,1.5*\j);
            \draw (1.5*\i,1.5*\j) -- (1.5*\i+1.5, 1.5*\j+1.5);
            \draw (1.5*\i+1.5,1.5*\j) -- (1.5*\i, 1.5*\j+1.5);
        }
    }
    \draw[color=blue, thick] (1.2, 1.5) -- (13.386, 13.677);
    \draw[color=red, thick] (1.3, 10.993) -- (11, 1.30015);
\end{tikzpicture}

    \caption{A chessboard with all of its triangles labeled. This chessboard is used to interpret the $n=8$ problem.
    The red line represents one of the diagonals of the chessboard.
    The blue line represents one of the anti-diagonals of the chessboard.
    }\label{f-chessboard}
\end{figure}

\subsection{Common notation}
In this section we describe variables and notation that are shared by the
lower and upper bound problems.

\paragraph{Chessboard triangle variables.}
In both problems, the variable $x$ consists of
$4$ $n \times n$ matrices $N, E, S, W$, 
and some additional slack variables.
The $i,j$th entry in $N, E, S, W$ is interpreted as a value associated with 
the North, East, South, or West triangle, respectively, formed by dividing 
each square of an $n\times n$ chessboard into $4$ right triangles.
We index these matrices starting at $0$, 
diverging from the notation used in \cite{simkin2021number,knuth2021xqueens} 
which use indexing that begins at $1$. 
Figure \ref{f-chessboard} shows the $n=8$ case.

\paragraph{Diagonal sum operators.}
We introduce operators $\mathcal D_k : \reals^{n\times n} \to \reals$ 
and $\mathcal A_k : \reals^{n\times n} \to \reals$
defined for $k \in \{-n, -n+1, -n +2, \ldots, -1, 0, 1, \ldots, n-2, n-1, n\}$,
where $\mathcal D_k(Z)$ is the sum of the $k$th diagonal of $Z$, and
$\mathcal A_k(Z)$ is the sum of the $k$th anti-diagonal of $Z$.
For example $\mathcal D_0 Z = \sum_{i=0}^{n-1} Z_{ii} = \tr{Z}$,
$\mathcal D_1 Z = \sum_{i=1}^{n-1} Z_{i,i-1}$, and
$\mathcal A_{-1} Z = \sum_{i=1}^{n-1} Z_{n-i,i}$.
Note that $\mathcal D_{n} Z =\mathcal D_{-n} Z = 0$.
These are illustrated in figure \ref{f-chessboard}:
$\mathcal D_{-2} N$ is the sum of the entries in the North triangles 
the red line passes through, and
$\mathcal A_1 E$ is the sum of the entries in the East triangles 
the blue line passes through.

\paragraph{Negative entropy.}
Following \cite{knuth2021xqueens}, we define the function $g:\reals_+\to \reals$
as $g(x)=x \log x$ for $x>0$, and $g(0) = 0$.  (This is the negative entropy
function \cite[p.72]{bv2004convex}.)

\subsection{Lower bound problem}\label{s-lower-bound-problem}

This problem formulation is taken from \cite[Claim 6.3]{simkin2021number},
except that Simkin maximizes a concave function and we minimize its negative, a
convex function.

\paragraph{Slack variables.}
We introduce the following slack variables,
\[
    d_k = \frac{1}{n} - \mathcal D_{k} (S+W) - \mathcal D_{k+1} (N+E),
    \quad k \in \{-n, -n+1, \ldots, n - 1\},
\]
and
\[
    a_k = \frac{1}{n} - \mathcal A_{k} (S+E) - \mathcal A_{k+1} (N+W),
    \quad k \in \{-n, -n+1, \ldots, n - 1\}.
\]
The quantities $d_k$ and $a_k$ are the sums along the diagonals and anti-diagonals,
respectively, of the chessboard.
In figure \ref{f-chessboard}, $d_{-1}$ includes contributions from all
triangles the red line passes through and $a_{0}$ includes contributions 
from triangles the blue line passes through.

These equations form $4n$ entries in $A, b$.

\paragraph{Objective.}
The objective function is \[
    \sum_{i=0}^{n-1} \sum_{j=0}^{n-1}
        \left(g(N_{i,j}) + g(E_{i,j}) + g(S_{i,j}) + g(W_{i,j})\right)
    + 
    \sum_{k=-n}^{n-1}
        \left(g(d_k) + g(a_k)\right)
    +
    4 \log n + 2 \log 2 + 3.
\]

\paragraph{Constraints.}
Simkin also introduces the constraints
\[
    \sum_{j=0}^{n-1} N_{i,j}+E_{i,j}+S_{i,j}+W_{i,j} = \frac{1}{n},
    \quad i \in \{0, 1, \ldots, n-1\},
\]
and \[
    \sum_{i=0}^{n-1} N_{i,j}+E_{i,j}+S_{i,j}+W_{i,j} = \frac{1}{n},
    \quad j \in \{0, 1, \ldots, n-1\}.
\]

These constraints are linearly dependent, with co-rank one,
so we delete the first constraint with
$i = 0$ to obtain a total of $2n - 1$ constraints that we include in $A, b$.

\paragraph{Properties.}
This problem has $n^2$ entries in each of $N, E, S, W$ and $2n$ entries in each
of $d, a$.
Accordingly, the total number of variables is $p=4n^2 + 4n$.
We have $4n$ constraints affecting the slack variables, and $2n - 1$ constraints
affecting only $N, E, S, W$ for a total of $q=  6n - 1$ constraints.

The objective is a sum of the negative entropy of individual optimization variables,
making it separable and strictly convex.

The constraint with the most variables are the row and column constraints, which
involve $4n$ variables.
Each triangle is in at most $1$ column constraint, $1$ row constraint,
$1$ diagonal constraint, and $1$ anti-diagonal constraint.
Therefore, each column of $A$ can have at most $4$ entries.

\subsection{Upper bound problem}

We use Knuth's formulation of the Xqueenon problem \cite{knuth2021xqueens},
except that he maximizes a concave function and we minimize its negative, a 
convex function.

\paragraph{Slack variables.}
We introduce the slack variables 
\[
d^{SW}_{k} = 1 -\frac{1}{2n} \mathcal D_{k}\left(S + W\right), \quad
k \in \{-n+1, -n+2, \ldots, n-1\},
\]
\[
d^{NE}_{k} = 1 -\frac{1}{2n} \mathcal D_{k}\left(N + E\right), \quad 
k \in \{-n+1, -n+2, \ldots, n -1\},
\]
\[
a^{SE}_{k} = 1 -\frac{1}{2n} \mathcal A_{k}\left(S + E\right), \quad
k \in \{-n+1, -n+2, \ldots, n-1\},
\]
and
\[
a^{NW}_{k} = 1 -\frac{1}{2n} \mathcal A_{k}\left(N + W\right), \quad 
k \in \{-n+1, -n+2, \ldots, n -1\}.
\]
These equations form $8n - 4$ entries in $A, b$.

\paragraph{Objective.}
For ease of notation, let
\[
    d^{SW}_{-n} = d^{NE}_{n} = a^{SE}_{-n} = a^{NW}_{n} = 1.
\]
Our objective function is
\[
        3 + L_0(N, E, S, W) + L_-(d^{SW}, d^{NE}) + L_+(a^{SE}, a^{NW}),
\]
where
\[
    L_0(N, E, S, W) = \frac{1}{4n^2} \sum_{i=0}^{n-1}\sum_{j=0}^{n-1}
        \left(g(N_{i,j}) + g(E_{i,j}) +g(S_{i,j}) +g(W_{i,j})\right),
\]
\[
    L_-(d^{SW}, d^{NE}) =
    \frac{1}{n}\sum_{k=-n+1}^{n}
        \int_{0}^1
            g\left((1 - y)d_{k-1}^{SW} + y d_{k}^{NE}\right)
        \,dy,
\]
and
\[
    L_+(a^{SE}, a^{NW}) =
    \frac{1}{n}\sum_{k=-n+1}^{n}
        \int_{0}^1
            g\left((1 - y)a_{k-1}^{SE} + y a_{k}^{NW}\right)
        \,dy.
\]

Using a symbolic solver, we were able to generate closed-form expressions for the
integrals and their partial derivatives \cite{sagemath}.

In order to make the matrix block-diagonal, $d^{SW}_k$ and $d^{NE}_k$ must
be interleaved in $x$.
Similar interleaving applies to $a^{SE}_k$ and $a^{NW}_k$.

\paragraph{Constraints}
In addition to the $8n - 4$ equations involving the slack variables, Knuth
requires the following conditions on $N$ and $S$,
\[
    \sum_{j=0}^{n-1} N_{i,j} = n, \quad i \in \{0, 1, \dots, n-1\},
\]
\[
    \sum_{j=0}^{n-1} S_{i,j} = n, \quad i \in \{0, 1, \dots, n-1\},
\]
and
\[
    \sum_{i=0}^{n-1} N_{i,j} + S_{i,j} = 2n, \quad j \in \{0, 1, \dots, n-1\}.
\]
As any of these equations are linearly dependent on all the others, we choose to
eliminate the first column constraint on $N$.

On $E$ and $W$, Knuth requires
\[
    \sum_{i=0}^{n-1} E_{i,j} = n, \quad j \in \{0, 1, \dots, n-1\},
\]
\[
    \sum_{i=0}^{n-1} W_{i,j} = n, \quad j \in \{0, 1, \dots, n-1\},
\]
and 
\[
    \sum_{j=0}^{n-1} E_{i,j} + W_{i,j} = 2n, \quad i \in \{0, 1, \dots, n-1\}.
\]
As with $N$ and $S$, one of these equations is linearly dependent, and we choose 
to eliminate the first row constraint on $E$.

\paragraph{Properties.}
This problem has $n^2$ entries in each of $N, E, S, W$ and $2n- 1$ entries in
each of $d^{SW}_k, d^{NE}_k, a^{SE}_k, a^{NW}_k$.
This forms a total of $p = 4n^2 + 8n - 4$ variables.
We also have the $8n - 4$ constraints involving the slack variables and $6n - 2$
of the other constraints for a total of $q = 14n - 6$ constraints.
The objective is block separable as each variable appears in only one term
of the objective function and no term has more than two variables.
The rows of $A$ with the most entries are the entries along the diagonal and
anti-diagonal, which contain $2n$ entries of $N,E,S,W$ and $1$ slack variable.
Each column of $A$ has at most $4$ non-zero entries: $1$ from its row constraint,
$1$ from its column constraints, $1$ from its diagonal term, and $1$ from its
anti-diagonal term.
Columns associated with slack variables have one non-zero entry.

\subsection{Approximate upper bound problem}\label{e-approx-upper-bnd}

In the initial phase of computing $U_n$ we  solve a problem with a diagonal
Hessian that approximates the upper bound problem.
We do this by applying Jensen's inequality to the integrals in the diagonal and
anti-diagonal terms of the objective.

After applying this approximation, the integral terms of the objective function become
\[
    g\left(\frac{1}{2}d^{SW}_{k-1} +\frac{1}{2}d^{NE}_k\right),
\]
and
\[
    g\left(\frac{1}{2}a^{SE}_{k-1} +\frac{1}{2}a^{NW}_k\right).
\]
We introduce new slack variables 
$d_k = \frac{1}{2}d^{SW}_{k-1} +\frac{1}{2}d^{NE}_k$ and
$a_k = \frac{1}{2}a^{SE}_{k-1} +\frac{1}{2}a^{NW}_k$ and then replace the
integral terms with $g(d_k)$ and $g(a_k)$ appropriately.

All other constraints and terms of the objective function are the same.

\end{document}